\documentclass[11pt]{article}
\usepackage{amssymb,amsfonts,amsmath, psfrag}
\usepackage{graphicx, colordvi, cite}
\usepackage{mathrsfs}
\usepackage{rotating}
\usepackage{hyperref}
\newtheorem{theo}{Theorem}

\newtheorem{defi}[theo]{Definition}
\newtheorem{lem}[theo]{Lemma}
\newtheorem{conj}[theo]{Conjecture}
\newtheorem{ex}{Example}

\makeatletter \@addtoreset{equation}{section}
\@addtoreset{theo}{section} \makeatother

\def\qed{\hfill \rule{4pt}{7pt}}
\def\pf{\noindent {\it Proof.} }

\def\lc{\mathop{\rm lc}}

\def\Sab{S_{a,b}}

\begin{document}

\begin{center}
{\large \bf
Polynomial Reduction and Super Congruences}

\vspace{0.3cm}

Qing-Hu Hou\\[8pt]

{\tt qh\_hou@tju.edu.cn} \\
School of Mathematics\\
Tianjin University\\
Tianjin 300350, P. R. China\\[12pt]

Yan-Ping Mu\\[8pt]
{\tt yanping.mu@gmail.com} \\
College of Science\\
Tianjin University of Technology\\
Tianjin 300384, P. R. China\\[6pt]

Doron Zeilberger \\[8pt]
{\tt doronzeil@gmail.com} \\
Department of Mathematics\\
Rutgers University\\
Piscataway, NJ 08854, USA\\[6pt]

\end{center}

\begin{abstract}
Based on a reduction processing, we rewrite a hypergeometric term as the sum of the difference of a hypergeometric term and a reduced hypergeometric term (the reduced part, in short). We show that when the initial hypergeometric term has a certain kind of symmetry, the reduced part contains only odd or even powers. As applications, we derived two infinite families of super-congruences.
\end{abstract}

\section{Introduction}
In recent years, many super congruences involving combinatorial sequences are discovered, see for example, Sun \cite{Sun12a}.
The standard methods for proving these congruences include combinatorial identities \cite{Sun14}, Gauss sums \cite{Alg}, symbolic computation \cite{Osb} et al.

We are interested in the following super congruence conjectured by van Hamme \cite{vH}
\[
\sum_{k=0}^{\frac{p-1}{2}} (-1)^k (4k+1) \left( \frac{(1/2)_k}{(1)_k} \right)^3 \equiv (-1)^{\frac{p-1}{2}} p \pmod{p^3},
\]
where $p$ is a odd prime and $(a)_k = a(a+1) \cdots (a+k-1)$ is the rising factorial.
This congruence was proved by Mortenson \cite{Mor}밃 Zudilin \cite{Zud} and Long \cite{Long11} by different methods.
Sun \cite{Sun12} proved a stronger version for prime $p \ge 5$
\[
\sum_{k=0}^{\frac{p-1}{2}} (-1)^k (4k+1) \left( \frac{(1/2)_k}{(1)_k} \right)^3 \equiv (-1)^{\frac{p-1}{2}} p + p^3 E_{p-3} \pmod{p^4},
\]
where $E_n$ is the $n$-th Euler number defined by
\[
\frac{2}{e^x+e^{-x}} = \sum_{n=0}^\infty E_n \frac{x^n}{n!}.
\]
A similar congruence was given by van Hamme \cite{vH} for $p \equiv 1 \pmod{4}$:
\[
\sum_{k=0}^{\frac{p-1}{2}} (4k+1) \left( \frac{(1/2)_k}{(1)_k} \right)^4 \equiv  p \pmod{p^3}.
\]
Long \cite{Long11} showed that in fact the above congruence holds for arbitrary odd prime modulo $p^4$.
Motivated by these two congruences, Guo \cite{Guo} proposed the following conjectures (corrected version).
\begin{conj}
\begin{itemize}
\item
For any odd prime $p$, positive integer $r$ and odd integer $m$, there exists an integer $a_{m.p}$ such that
\begin{equation}\label{conj3}
\sum_{k=0}^{\frac{p^r-1}{2}} (-1)^k (4k+1)^m \left( \frac{(1/2)_k}{(1)_k} \right)^3 \equiv a_{m,p} p^r (-1)^{\frac{(p-1)r}{2}}  \pmod{p^{r+2}}.
\end{equation}

\item
For any odd prime $p > (m-1)/2$, positive integer $r$ and odd integer $m$, there exists an integer $b_{m,p}$ such that
\begin{equation}\label{conj4}
\sum_{k=0}^{\frac{p^r-1}{2}} (4k+1)^m \left( \frac{(1/2)_k}{(1)_k} \right)^4 \equiv b_{m,p} p^r  \pmod{p^{r+3}}.
\end{equation}
\end{itemize}
\end{conj}
Liu \cite{Liu} and Wang \cite{Wang} confirmed the conjectures for $r=1$ and some initial values $m$.
Jana and Kalita \cite{Jana} and Guo \cite{Guo2} confirmed \eqref{conj3} for $m = 3$ and $r \ge 1$.
We will prove a stronger version of \eqref{conj3} for the case of $r=1$ and arbitrary odd $m$ and a weaker version of \eqref{conj4} for the case of $r=1$ and arbitrary odd $m$ by a reduction process.

Recall that a hypergeometric term $t_k$ is a function of $k$ such that $t_{k+1}/t_k$ is a rational function of $k$.
Our basic idea is to rewrite the product of a polynomial $f(k)$ in $k$ and a hypergeometric term $t_k$ as
\[
f(k) t_k = \Delta ( g(k)t_k ) + h(k)t_k = (g(k+1)t_{k+1}-g(k)t_k) +  h(k)t_k,
\]
where $g(k), h(k)$ are polynomials in $k$ such that the degree of $h(k)$ is bounded. To this aim, we construct $x(k)$ such that $\Delta x(k)t_k$ equals the product of $t_k$ and a polynomial $u(k)$ and that $f(k)$ and $u(k)$ has the same leading term. Then we have
\[
f(k)t_k -  \Delta x(k)t_k = (f(k)-u(k)) t_k
\]
is the product of $t_k$ and a polynomial of degree less than $f(k)$. We call such a reduction process one reduction step. Continuing this reduction process, we finally obtain a polynomial $h(k)$ with bounded degree. We will show that for $t_k=\left( \frac{(1/2)_k}{(1)_k} \right)^r$, $r=3,4$ and an arbitrary polynomial of form $(4k+1)^m$ with $m$ odd, the reduced polynomial $h(k)$ can be taken as $(4k+1)$. This enables us to reduce the congruences \eqref{conj3} and \eqref{conj4} to the special case of $m=1$, which is known for $r=1$.

We notice that  Pirastu-Strehl \cite{Pir} and Abramov \cite{Abr75, Abr95} gave the minimal decomposition when $t_k$ is a rational function, Abramov-Petkov\v{s}ek \cite{Abr01, Abr02} gave the minimal decomposition when $t_k$ is a hypergeometric term, and Chen-Huang-Kauers-Li \cite{Chen} applied the reduction to give an efficient creative telescoping algorithm. These algorithms concern a general hypergeometric term. While we focus on a kind of special hypergeometric term so that the reduced part $h(k)t_k$ has a nice form.

The paper is organized as follows. In Section 2, we consider the reduction process for a general hypergeometric term $t_k$. Then in Section 3 we consider those $t_k$ with the property $a(k)$ is a shift of $-b(k)$, where $t_{k+1}/t_k=a(k)/b(k)$. As an application, we prove a stronger version of \eqref{conj3} for the case $r=1$. Finally, we consider the case of $a(k)$ is a shift of $b(k)$, which corresponds to \eqref{conj4}. In this case, we show
that there is a rational number $b_m$ instead of an integer such that \eqref{conj4} holds when $r=1$.

\section{The Difference Space and Polynomial Reduction}

Let $K$ be a field and $K[k]$ be the ring of polynomials in $k$ with coefficients in $K$.
Let $t_k$ be a hypergeometric term. Suppose that
\[
\frac{t_{k+1}}{t_k} =  \frac{a(k)}{b(k)},
\]
where $a(k), b(k)\in K[k]$. It is straightforward to verify that
\begin{equation}\label{Gosper}
\Delta_k \left( b(k-1) x(k) t_k \right) = (a(k) x(k+1) - b(k-1) x(k)) t_k.
\end{equation}
We thus define the {\it difference space} corresponding to $a(k)$ and $b(k)$ to be
\[
\Sab = \{a(k) x(k+1) - b(k-1) x(k) \colon x(k) \in K[k] \}.
\]
We see that for $f(k) \in \Sab$, we have $f(k)t_k = \Delta_k(p(k)t_k)$ for a  certain polynomial $p(k) \in K[k]$.

Let $\mathbb{N}, \mathbb{Z}$ denote the set of nonnegative integers and the set of integers, respectively.
Given $a(k), b(k) \in K[k]$, we denote
\begin{equation}\label{u}
u(k) = a(k) - b(k-1),
\end{equation}
\begin{equation}\label{d}
d = \max\{ \deg u(k), \deg a(k)-1 \},
\end{equation}
and
\begin{equation}\label{m0}
m_0 = - \lc u(k) / \lc a(k),
\end{equation}
where $\lc p(k)$ denotes the leading coefficient of $p(k)$.

We first introduce the concept of degeneration.

\begin{defi}
Let $a(k),b(k) \in K[k]$ and $u(k), m_0$ be given by \eqref{u} and \eqref{m0}.
If
\[
\deg u(k) = \deg a(k) - 1 \quad \mbox{and} \quad  m_0 \in \mathbb{N},
\]
we say that the pair $(a(k),b(k))$ is \emph{degenerated}.
\end{defi}

We will see that the degeneration is closely related to the degrees of the elements in $\Sab$.
\begin{lem}\label{deg}
Let $a(k),b(k) \in K[k]$ and $d, m_0$ be given by \eqref{d} and \eqref{m0}.  For any polynomial $x(k) \in K[k]$, let
\[
p(k) = a(k) x(k+1) - b(k-1) x(k).
\]
If $(a(k),b(k))$ is degenerated and $\deg x(k)=m_0$, then $\deg p(k)<d+m_0$; Otherwise, $\deg p(k) = d+\deg x(k)$.
\end{lem}
\pf
Notice that
\[
p(k) = u(k) x(k) + a(k) (x(k+1)-x(k)).
\]
If the leading terms of $u(k) x(k)$ and $a(k) (x(k+1)-x(k))$ do not cancel, the degree of $p(k)$ is $d + \deg x(k)$. Otherwise, we have $\deg u(k) = \deg a(k)-1$ and
\[
\lc u(k) + \lc a(k) \cdot \deg x(k) = 0,
\]
i.e., $\deg x(k) = m_0$. \qed

It is clear that $\Sab$ is a subspace of $K[k]$, but is not a sub-ring of $K[k]$ in general.
Let $[p(k)]=p(k)+\Sab$ denote the coset of a polynomial $p(k)$.
We see that the quotient space $K[k]/\Sab$ is  finite dimensional.

\begin{theo}
Let $a(k),b(k) \in K[k]$ and $d, m_0$ be given by \eqref{d} and \eqref{m0}. We have
\[
K[k]/\Sab
= \begin{cases}
  \langle [k^0], [k^1], \ldots, [k^{d-1}], [k^{d+m_0}] \rangle, & \mbox{if $(a(k),b(k))$ is degenerated,} \\[5pt]
 \langle [k^0], [k^1], \ldots, [k^{d-1}] \rangle, & \mbox{otherwise}.
\end{cases}
\]
\end{theo}

\pf
For any nonnegative integer $s$, let
\[
p_s(k) = a(k) (k+1)^s - b(k-1) k^s.
\]

We first consider the case when the pair $(a(k),b(k))$ is not degenerated.
By Lemma~\ref{deg}, we have
\[
\deg p_s(k) = d + s, \quad \forall\, s \ge 0.
\]
Suppose that $p(k)$ is a polynomial of degree $m \ge d$. Then
\begin{equation}\label{reduction}
p'(k) = p(k) - \frac{\lc p(k)}{\lc p_{m-d}(k)} p_{m-d}(k)
\end{equation}
is a polynomial of degree less than $m$ and $p(k) \in [p'(k)]$.
By induction on $m$, we derive that for any polynomial $p(k)$ of degree $\ge d$, there exists a polynomial $\tilde{p}(k)$ of degree $<d$ such that $p(k) \in [\tilde{p}(k)]$. Therefore,
\[
K[k]/\Sab = \langle [k^0], [k^1], \ldots, [k^{d-1}] \rangle.
\]

Now assume that $(a(k),b(k))$ is degenerated. By Lemma~\ref{deg},
\[
\deg p_s(k) = d+s, \quad \forall\, s \not= m_0 \quad \mbox{and} \quad
\deg p_{m_0}(k) < d+m_0.
\]
The above reduction process \eqref{reduction} works well except for the polynomials $p(k)$ of degree $d+m_0$. But in this case,
\[
p(k) - \lc p(k)  \cdot k^{d+m_0}
\]
is a polynomial of degree less than $d+m_0$. Then the reduction process continues until the degree is less than $d$. We thus derive that
\[
K[k]/\Sab = \langle [k^0], [k^1], \ldots, [k^{d-1}], [k^{d+m_0}] \rangle,
\]
completing the proof. \qed

\begin{ex}
Let $n$ be a positive integer and
\[
t_k = (-n)_k/k!,
\]
where $(\alpha)_k = \alpha(\alpha+1) \cdots (\alpha+k-1)$ is the raising factorial. Then
\[
a(k) = k-n, \quad b(k)=k+1,
\]
and
\[
\Sab = \{ (k-n) \cdot x(k+1) - k \cdot x(k) \colon x(k) \in K[k] \}.
\]
We have
\[
K[k]/ \Sab = \langle [k^n] \rangle
\]
is of dimension one.
\end{ex}

\section{The case when $a(k)=-b(k+\alpha)$}
In this section, we consider the case when $a(k)=-b(k+\alpha)$ and $b(k)$ has a symmetric property. We will show that in this case, the reduction process maintains the symmetric property. Notice that in this case
\[
u(k) = a(k)-b(k-1)=-b(k+\alpha)-b(k-1)
\]
has the same degree as $a(k)$, the pair $(a(k),b(k))$ is not degenerated.

We first consider the relation between the symmetric property and the expansion of a polynomial.

\begin{lem}\label{sym}
Let $p(k) \in K[k]$ and $\beta \in K$. Then the following two statements are equivalent.
\begin{itemize}
\item[{\rm (1)}]
$p(\beta+k) = p(\beta-k)$ ($p(\beta+k)=-p(\beta-k)$, respectively).

\item[{\rm (2)}]
$p(k)$ is the linear combination of $(k-\beta)^{2i},\, i=0,1,\ldots$ ($(k-\beta)^{2i+1},\, i=0,1,\ldots$, respectively).
\end{itemize}
\end{lem}
\pf
Suppose that
\[
p(\beta+k) = \sum_i c_i k^i.
\]
Then
\[
p(\beta-k) = \sum_i c_i (- k)^i.
\]
Therefore,
\[
p(\beta+k)=p(\beta-k) \Longleftrightarrow c_{2i+1}=0,\ i=0,1,\ldots.
\]
The case of $p(\beta+k)= - p(\beta-k)$ can be proved in a similar way. \qed

Now we are ready to state the main theorem.
\begin{theo}\label{a+b}
Let $a(k),b(k) \in K[k]$ such that
\[
a(k) = - b(k+\alpha) \quad \mbox{and} \quad b(\beta+k) = \pm b(\beta-k),
\]
for some $\alpha,\beta \in K$.
Then for any non-negative integer $m$, we have
\[
[(k + \gamma)^{2m}] \in \left\langle [(k+ \gamma)^{2i}] \colon 0 \le 2 i < \deg a(k) \right\rangle
\]
and
\[
[(k + \gamma)^{2m+1}] \in \left\langle [(k+ \gamma)^{2i+1}] \colon 0 \le 2i+1 < \deg a(k) \right\rangle,
\]
where
\begin{equation}\label{ga}
\gamma = - \beta + \frac{\alpha-1}{2}.
\end{equation}
\end{theo}
\pf We only prove the case of $b(\beta+k)=b(\beta-k)$. The case of $b(\beta+k)=-b(\beta-k)$ can be proved in a similar way.
By Lemma~\ref{sym}, we may assume that
\[
b(k) = b_r (k-\beta)^r + b_{r-2} (k-\beta)^{r-2} + \cdots + b_0,
\]
where $r = \deg a(k) = \deg b(k)$ is even and $b_r, b_{r-2}, \ldots, b_0 \in K$ are the coefficients.

Since $(a(k),b(k))$ is not degenerated, taking
\begin{equation}\label{xs}
x(k) = x_s(k)= - \frac{1}{2} \left( k+ \gamma - \frac{1}{2} \right)^s
\end{equation}
in Lemma~\ref{deg}, we derive that
\begin{equation}\label{ps}
p_s(k) = a(k) x_s(k+1) - b(k-1) x_s(k)
\end{equation}
is a polynomial of degree $s+r$. More explicitly, we have
\[
p_s(k) =  \frac{1}{2} \left(b(k+\alpha) \left( k+ \gamma + \frac{1}{2} \right)^s
+ b(k-1) \left( k+ \gamma - \frac{1}{2} \right)^s \right)
\]
is a polynomial with leading term $b_r k^{s+r}$.

Notice that
\[
p_s(-\gamma+k) =  \frac{1}{2} \left( b(k+\alpha-\gamma) \left( k + \frac{1}{2} \right)^s
+ b(k-\gamma-1) \left( k - \frac{1}{2} \right)^s \right)
\]
and
\begin{align*}
p_s(-\gamma-k) & = \frac{1}{2}  \left( b(-k+\alpha-\gamma) \left( - k + \frac{1}{2} \right)^s
+ b(-k-\gamma-1) \left( - k - \frac{1}{2} \right)^s \right) \\
& =
\frac{(-1)^{s}}{2}  \left( b(-k+\alpha-\gamma) \left( k - \frac{1}{2} \right)^s
+ b(-k-\gamma-1) \left(  k + \frac{1}{2} \right)^s \right).
\end{align*}
Since $b(\beta+k)=b(\beta-k)$, i.e., $b(k)=b(2\beta-k)$, we deduce that
\begin{align*}
& p_s(-\gamma-k) \\
& = \frac{(-1)^{s}}{2} \left( b(2 \beta + k-\alpha+\gamma) \left( k - \frac{1}{2} \right)^s
+ b(2 \beta + k+ \gamma+1) \left(k + \frac{1}{2} \right)^s \right).
\end{align*}
By the relation \eqref{ga}, we derive that
\[
p_s(-\gamma-k) = (-1)^s p_s(-\gamma+k).
\]

Suppose that $p(k)$ is a linear combination of the even powers of $(k+\gamma)$ and $\deg p(k) \ge r$. By Lemma~\ref{sym}, we have $p(-\gamma-k)=p(-\gamma+k)$ and thus
\[
p'(k) = p(k) - \frac{\lc p(k)}{b_r} \cdot p_{\deg p(k) - r}(k)
\]
also satisfies $p'(-\gamma-k)=p'(-\gamma+k)$ since $\deg p(k)$ and $r$ are both even. It is clear that
$p(k) \in [p'(k)]$ and the degree of $p'(k)$ is less than the degree of $p(k)$.
Continuing this reduction process, we finally derive that $p(k) \in [\tilde{p}(k)]$ for some polynomial $\tilde{p}(k)$ with degree $< r$ and  satisfying $\tilde{p}(-\gamma-k)=\tilde{p}(-\gamma+k)$. Therefore,
\[
[p(k)] \in \langle [(k+\gamma)^{2i}] \colon 0 \le 2i < r \rangle.
\]

Suppose that $p(k)$ is a linear combination of the odd powers of $(k+\gamma)$ and $\deg p(k) \ge r$. Then we have
$p(-\gamma-k)=-p(-\gamma+k)$ and thus
\[
p'(k) = p(k) - \frac{\lc p(k)}{b_r} \cdot p_{\deg p(k) - r}(k)
\]
also satisfies $p'(-\gamma-k)=-p'(-\gamma+k)$. Continuing this reduction process, we finally derive that
\[
[p(k)] \in \langle [(k+\gamma)^{2i+1}] \colon 0 \le 2i+1 < r \rangle.
\]
This completes the proof. \qed

We may further require to express $[(k+ \gamma)^m]$ as an integral linear combination of $[(k+ \gamma)^i],\ 0 \le i <r$ when $b(k)=(k+1)^r$.
\begin{theo}\label{int-a+b}
Let
\[
t_k = (-1)^k \left( \frac{(\alpha)_k}{k!} \right)^{r},
\]
where $r$ is a positive integer and $\alpha$ is a rational number with denominator $D$.
Then for any positive integer $m$, there exist integers $a_0,\ldots,a_{r-1}$ and a polynomial $x(k) \in \mathbb{Z}[k]$ such that
\[
(2Dk+D \alpha)^m t_k = \sum_{i=0}^{r-1} a_i (2 D k+ D \alpha)^i t_k + \Delta_k \left( 2^{r-1} (Dk)^{r} x(2Dk) t_k \right).
\]
Moreover, $a_i=0$ if $i \not\equiv m \pmod{2}$.
\end{theo}
\pf
We have
\[
\frac{t_{k+1}}{t_k} = \frac{-(k+\alpha)^r}{(k+1)^r}.
\]
Let
\[
a(k)=-(k+\alpha)^r \quad \mbox{and} \quad b(k)=(k+1)^r.
\]
We see that it is the case of $\beta=-1$ and $\gamma=\alpha/2$ of Theorem~\ref{a+b}.
From \eqref{Gosper}, we derive that
\begin{equation}\label{Delta}
\Delta_k( k^r x_s(k) t_k) = p_s(k) t_k,
\end{equation}
where $x_s(k)$ and $p_s(k)$ are given by \eqref{xs} and \eqref{ps} respectively.
Multiplying $(2D)^{s+r}$ on both sides, we obtain
\begin{equation}\label{xs-ps}
\Delta_k(2^{r-1} (Dk)^r \tilde{x}_s(2Dk) t_k) = \tilde{p}_s(k') t_k,
\end{equation}
where $k'=2Dk+D\alpha$,
\begin{equation}\label{xsp}
\tilde{x}_s(k) = - (k+D \alpha - D)^s,
\end{equation}
and
\begin{equation}\label{psp}
\tilde{p}_s(k) = \frac{1}{2} \left( (k+ D\alpha)^r \left( k+D \right)^s + (k - D \alpha)^r \left( k-D \right)^s \right).
\end{equation}
Notice that $\tilde{x}_s(k), \tilde{p}_s(k) \in \mathbb{Z}[k]$ and $\tilde{p}_s(k)$ is a monic polynomial of degree $s+r$.
Moreover, $\tilde{p}_s(k)$ contains only even powers of $k$ or only odd powers of $k$. Using $\tilde{p}_s(k)$ to do the reduction \eqref{reduction}, we derive that there exists integers $c_m, c_{m-2}, \ldots$ such that
\[
p(k) = k^m - c_m \tilde{p}_{m-r}(k) - c_{m-2} \tilde{p}_{m-r-2}(k) - \cdots
\]
becomes a polynomial of degree less than $r$. Clearly, $p(k) \in \mathbb{Z}[k]$.
Replacing $k$ by $k'$ and multiplying $t_k$, we derive that
\[
(k')^m t_k = p(k') t_k + \Delta_k (2^{r-1} (Dk)^r (c_m \tilde{x}_{m-r}(2Dk) + c_{m-2} \tilde{x}_{m-r-2}(2Dk) + \cdots) t_k ),
\]
completing the proof. \qed

%
%

As an application, we confirm Conjecture 6 of \cite{Liu}.

\begin{theo}
Let
\[
S_m = \sum_{k=0}^{\frac{p-1}{2}} (-1)^k (4k+1)^m \left( \frac{(1/2)_k}{(1)_k} \right)^3.
\]
For any positive odd integer $m$, there exist integers $a_m$ and $c_m$ such that
\[
S_m \equiv a_m \left( p (-1)^{\frac{p-1}{2}} + p^3 E_{p-3} \right) + p^3 c_m \pmod{p^4}
\]
holds for any prime $p \ge 5$.
\end{theo}
\pf
Taking $r=3$ and $\alpha=1/2$ in Theorem~\ref{int-a+b}, there exist an integer $a_m$ and a polynomial $q_m(k) \in \mathbb{Z}[k]$ such that
\[
(4k+1)^m t_k  -  a_m (4k+1) t_k = \Delta_k ( 32 k^3 q_m(4k) t_k ) ,
\]
where $t_k = (-1)^k (\frac{1}{2})_k^3/(1)_k^3$.
Summing over $k$ from $0$ to $\frac{p-1}{2}$, we derive that
\[
S_m - a_m S_1 = 32 \omega^3 q_m(4 \omega) (-1)^\omega  \left( \frac{(1/2)_\omega}{(1)_\omega} \right)^3,
\]
where $\omega=\frac{p+1}{2}$.
Noting that
\[
\frac{(1/2)_\omega}{(1)_\omega} = p \frac{1}{p+1} \prod_{i=1}^{\frac{p-1}{2}} \frac{2i-1}{2i}
\]
and
\[
\frac{1}{p+1} \prod_{i=1}^{\frac{p-1}{2}} \frac{2i-1}{2i} = \frac{1}{p+1} \prod_{i=1}^{\frac{p-1}{2}} \frac{p-2i}{2i} \equiv (-1)^{\frac{p-1}{2}} \pmod{p},
\]
we have
\[
\left( \frac{(1/2)_\omega}{(1)_\omega} \right)^3 \equiv p^3 (-1)^{\frac{p-1}{2}} \pmod{p^4}.
\]
Hence
\[
S_m - a_m S_1 \equiv  - 32 p^3 \omega^3 q_m(4 \omega) \pmod{p^4}
\]
Let $c_m =  - 4 q_m(2)$. We then have
\[
S_m \equiv a_m S_1 + p^3 c_m \pmod{p^4}.
\]
Sun \cite{Sun12} proved that for any prime $p \ge 5$,
\[
S_1 \equiv (-1)^{\frac{p-1}{2}} p + p^3 E_{p-3} \pmod{p^4}.
\]
Therefore,
\[
S_m \equiv a_m \left( p (-1)^{\frac{p-1}{2}} + p^3 E_{p-3} \right) + p^3 c_m \pmod{p^4}. \tag*{\qed}
\]

\noindent {\it Remark 1.} The coefficient $a_m$ and the polynomial $q_m(k)$ can be computed by the extended Zeilberger's algorithm \cite{EZ}.

\section{The case when $a(k)=b(k+\alpha)$}

We first give a criterion on the degeneration of $(a(k),b(k))$.
\begin{lem}\label{degen}
Let $a(k),b(k) \in K[k]$ such that $a(k) = b(k+\alpha)$. Suppose that $-(\alpha+1) \deg a(k) \not\in \mathbb{N}$.
Then $(a(k),b(k))$ is not degenerated.
\end{lem}
\pf
Let $r = \deg a(k) = \deg b(k)$ and
\[
u(k) = a(k) - b(k-1) = b(k+\alpha) - b(k-1).
\]
It is clear that the coefficient of $k^r$ in $u(k)$ is $0$ and
the coefficient of $k^{r-1}$ in $u(k)$ is $\lc b(k) \cdot  (\alpha+1) r$. Since $ (\alpha+1) r \not= 0$, we derive that $\deg u(k) = r - 1$. Thus,
\[
- \lc u(k)/\lc a(k) = - \lc u(k)/\lc b(k) = - (\alpha+1) r.
\]
Since $-(\alpha+1)r \not\in \mathbb{N}$, the pair $(a(k),b(k))$ is not degenerated.
\qed

When $a(k)$ is a shift of $b(k)$, we have a result similar to Theorem~\ref{a+b}.
\begin{theo}\label{a-b}
Let $a(k),b(k) \in K[k]$ such that
\[
a(k) = b(k+\alpha) \quad \mbox{and} \quad b(\beta+k) = \pm b(\beta-k),
\]
for some $\alpha,\beta \in K$.
Assume further that $-(\alpha+1) \deg a(k) \not\in \mathbb{N}$. Then for any non-negative integer $m$, we have
\[
(k + \gamma)^{2m} \in \left\langle [(k+ \gamma)^{2i}] \colon 0 \le 2 i < \deg a(k) - 1  \right\rangle
\]
and
\[
(k + \gamma)^{2m+1} \in \left\langle [(k+ \gamma)^{2i+1}] \colon 0 \le 2i+1 < \deg a(k)-1 \right\rangle,
\]
where
\[
\gamma = - \beta + \frac{\alpha-1}{2}.
\]
\end{theo}
\pf The proof is parallel to the proof of Theorem~\ref{a+b}. Instead of \eqref{xs}, we take
\[
x(k) = x_s(k)= \left( k+ \gamma - \frac{1}{2} \right)^s
\]
in Lemma~\ref{deg}. By Lemma~\ref{degen}, $(a(k),b(k))$ is not degenerated and
\[
\deg (a(k)-b(k-1)) = \deg a(k) - 1.
\]
Hence the polynomial
\[
p_s(k) = a(k) x_s(k+1) - b(k-1) x_s(k)
\]
satisfies
\[
\deg p_s(k) = s + \deg a(k) - 1.
\]
Moreover, we have
\[
p_s(-\gamma-k) = (-1)^{s+1} p_s(-\gamma+k),
\]
so that the reduction process maintains the symmetric property.
Therefore, the reduction process continues until the degree is less than $\deg a(k) -1$. \qed

Similar to Theorem~\ref{int-a+b}, we have the following result.

\begin{theo}\label{int-a-b}
Let
\[
t_k = \left( \frac{(\alpha)_k}{k!} \right)^{r},
\]
where $r$ is a positive integer and $\alpha$ is a rational number with denominator $D$.
Suppose that $-\alpha r \not\in \mathbb{N}$.
Then for any positive integer $m$, there exist integers $a_0,\ldots,a_{r-2}$ and a polynomial $x(k) \in \mathbb{Z}[k]$ such that
\[
(2Dk+D \alpha)^m t_k = \frac{1}{C_m} \sum_{i=0}^{r-2} a_i (2 D k+ D \alpha)^i t_k + \frac{1}{C_m} \Delta_k \left( 2^{r-1}(Dk)^{r} x(2Dk) t_k \right),
\]
where
\[
C_m = \prod_{0 \le 2i \le m-r+1} ( (\alpha r + m -r+1-2i) \cdot D) .
\]
Moreover, $a_i=0$ if $i \not\equiv m \pmod{2}$.
\end{theo}
\pf
The proof is parallel to the proof of Theorem~\ref{int-a+b}. Instead of \eqref{xsp} and \eqref{psp}, we take
\begin{equation}\label{xs1}
\tilde{x}_s(k) = \left( k + D \alpha - D \right)^s
\end{equation}
and
\begin{equation}\label{ps1}
\tilde{p}_s(k) =  \frac{1}{2} ( (k+ D\alpha)^r \left( k+D \right)^s - (k - D \alpha)^r \left( k-D \right)^s ),
\end{equation}
so that \eqref{xs-ps} still holds. It is clear that $\tilde{x}_s(k), \tilde{p}_s(k) \in \mathbb{Z}[k]$. But in this case, $\tilde{p}_s(k)$ is not monic. The leading term of $\tilde{p}_s(k)$ is
\[
(\alpha r + s) D \cdot k^{s+r-1}.
\]

Now let us consider the reduction process. Let $p(k) \in \mathbb{Z}[k]$ be a polynomial of degree $\ell \ge r-1$. Assume further that $p(k)$ contains only even powers of $k$ or only odd powers of $k$.
Setting
\begin{align*}
p'(k) & = \lc \tilde{p}_{\ell - r+1} (k) \cdot p(k) -\lc p(k) \cdot \tilde{p}_{\ell - r+1} (k) \\
& = (\alpha r + \ell-r+1) D \cdot p(k) -\lc p(k) \cdot \tilde{p}_{\ell - r+1} (k),
\end{align*}
we see that $p'(k) \in \mathbb{Z}[k]$ and $\deg p'(k)<\ell$. Since $\tilde{p}_{\ell-r+1}(k)$ contains only even powers of $k$ or only odd powers of $k$, so does $p'(k)$. Therefore, $\deg p'(k) \le \ell-2$.

Continuing this reduction process until $\ell < r-1$, we finally obtain that there exist integers $c_m, c_{m-2}, \ldots$ such that
\[
C_m k^m - c_m \tilde{p}_{m - r+1} (k) - c_{m-2} \tilde{p}_{m-r-1}(k) - \cdots,
\]
is a polynomial of degree less than $r-1$ and with integral coefficients,
where $C_m$ is the product of the leading coefficient of $\tilde{p}_{m- r+1} (k), \tilde{p}_{m-r-1}(k), \ldots$
\[
C_m = \prod_{0 \le 2i \le m-r+1} \left( (\alpha r + m-r+1-2i)D \right),
\]
as desired. \qed

For the special case of $t_k=(1/2)_k^4/(1)_k^4$, we may further reduce the factor $C_m$.
\begin{lem}\label{redu}
Let $m$ be a positive integer and
\[
t_k = \frac{(1/2)_k^4}{(1)_k^4}.
\]
\begin{itemize}
\item
If $m$ is odd, then there exist an integer $c$ and a polynomial $x(k) \in \mathbb{Z}[k]$  such that
\[
(4k+1)^m t_k = \frac{c}{C'_m} (4k+1) t_k + \frac{1}{C'_m} \Delta_k \left( 32 k^4 x(4k) t_k \right),
\]
where $C'_m = (\frac{m-1}{2})!$.

\item
If $m$ is even, then there exist integers $c,c'$ and a polynomial $x(k) \in \mathbb{Z}[k]$ such that
\[
(4k+1)^m t_k = \frac{1}{C'_m} (c + (4k+1)^2 c') t_k + \frac{1}{C'_m} \Delta_k \left( 64 k^4 x(4k) t_k \right),
\]
where $C'_m = (m-1)!!$.
\end{itemize}
\end{lem}
\pf This is the special case of Theorem~\ref{int-a-b} for $\alpha=1/2$ and $r=4$. Therefore, $D=2$ and $\alpha r - r + 1=-1$.

We need only to show that the coefficients of $\tilde{p}_s(k)$ given by \eqref{ps1} is divisible by $2$ when $s$ is odd and is divisible by $4$ when $s$ is even. Then we may replace $\tilde{x}_s(k)$ given by \eqref{xs1} by $\tilde{x}_s(k)/4$ and $\tilde{x}_s(k)/2$ so that the leading coefficient of $\tilde{p}_s(k)$ is reduced. Correspondingly, the product $C_m$ of the leading coefficients becomes
\[
\prod_{0 \le 2i \le m-3} \frac{1}{2} \lc \tilde{p}_{m-3-2i}(k) =  \prod_{0 \le 2i \le m-3}  (m-1-2i) = (m-1)!!, \quad \mbox{$m$ even,}
\]
and
\[
\prod_{0 \le 2i \le m-3} \frac{1}{4} \lc \tilde{p}_{m-3-2i}(k) =  \prod_{0 \le 2i \le m-3}  \frac{m-1-2i}{2} = \left(\frac{m-1}{2} \right)!, \quad \mbox{$m$ odd}.
\]

Notice that
\[
\tilde{p}_s(k) =  \frac{1}{2} ( (k+ 1)^4 \left( k+2\right)^s - (k - 1)^4 \left( k-2 \right)^s ).
\]
The coefficient of $k^j$ is
\[
\frac{1-(-1)^{s-j}}{2} \sum_{0 \le \ell \le 4, \ 0 \le j-\ell \le s} {4 \choose \ell} {s \choose j-\ell} 2^{s-j+\ell}.
\]
If $j-\ell<s$, the corresponding summand is divisible by $2$. If $j-\ell=s$ and $\ell$ is even, then $(-1)^{s-j}=1$ and the coefficient is $0$. Otherwise, $\ell=1$ or $\ell=3$, and thus $4 \mid {4 \choose \ell}$. Therefore, the coefficient must be divisible by $2$.

Now consider the case of $s$ being even. If $j-\ell<s-1$, the corresponding summand is divisible by $4$. Otherwise $j-\ell=s$ or $j-\ell=s-1$. We have seen that if $j-\ell=s$, then the coefficient is divisible by $4$. If $j-\ell=s-1$. Then
\[
 {s \choose j-\ell}=s  \quad \mbox{and} \quad 2^{s-j+\ell}=2.
\]
Thus the summand is also divisible by $4$. \qed

\begin{ex}
Consider the case of $m=11$. We have
\[
(4k+1)^{11} t_k + 10515 (4k+1) t_k = \Delta_k (32 k^4 p(k) t_k)
\]
where
\[
p(k) = \frac{1}{5} (4k-1)^8 - \frac{249}{20} (4k-1)^6 + \frac{20207}{60} (4k-1)^4 -\frac{89909}{20} (4k-1)^2 + \frac{524029}{20}.
\]
\end{ex}

As an application, we obtain the following congruences.
\begin{theo}
Let $m$ be a positive odd integer and $\mu=(m-1)/2$.
Denote
\[
S_m = \sum_{k=0}^{\frac{p-1}{2}} (4k+1)^m \left( \frac{(1/2)_k}{(1)_k} \right)^4.
\]
Then there exists an integer $a_m$ such that for each prime $p>\mu$,
\[
 S_m \equiv \frac{a_m}{\mu!} p \pmod{p^4}.
\]
\end{theo}
\pf
By Lemma~\ref{redu}, there exist an integer $a_m$ and a polynomial $q_m(k) \in \mathbb{Z}[k]$ such that
\[
(4k+1)^m t_k - \frac{a_m}{\mu!} (4k+1) t_k
= \frac{1}{\mu!} \Delta_k \left( 32 k^4 q_m(k) t_k \right),
\]
where $t_k = \left( \frac{(1/2)_k}{(1)_k} \right)^4$.
Summing over $k$ from $0$ to $(p-1)/2$, we obtain
\[
S_m - \frac{a_m}{\mu!} S_1 =
32 \omega^4 \frac{q_m(4 \omega)}{\mu!} \left( \frac{(1/2)_\omega}{(1)_\omega} \right)^4,
\]
where $\omega=(p+1)/2$.
When $p>\mu$, $1/\mu!$ is a $p$-adic integer and
\[
\left( \frac{(1/2)_\omega}{(1)_\omega} \right)^4 \equiv 0 \pmod{p^4}.
\]
Therefore,
\[
S_m \equiv \frac{a_m}{\mu!} S_1 \pmod{p^4}.
\]
It is shown by Long \cite{Long11} that
\[
S_1 \equiv p \pmod{p^4},
\]
completing the proof.  \qed


The integer $a_m$ and the polynomial $q_m(k)$ can be computed by the extended Zeilberger's algorithm.

By checking the initial values, we propose the following conjecture.
\begin{conj}
For any positive odd integer $m$, the coefficient $a_m/(\frac{m-1}{2})!$ is an integer.
\end{conj}

\vskip 15pt

{\small \noindent
{\bf Acknowledgement.} The work was supported by the National Natural Science Foundation of China (grants 11471244, 11771330 and 11701420).
}

\end{document}